\documentclass[a4paper,11pt]{article}
\usepackage{amsfonts,amssymb,latexsym,amsmath}
\usepackage{multicol}
\newcounter{num}[section]
\newcommand{\Num}{\refstepcounter{num}%
\textbf{\arabic{section}.\arabic{num}}}

\newcommand{\Theorem}{\textbf{Theorem~}}
\newcommand{\Proof}{\textbf{Proof}}
\newcommand{\Def}{\textbf{Definition~}}

\newcommand{\Lemma}{\textbf{Lemma~}}

\newcommand{\Prop}{\textbf{Proposition~}}
\newcommand{\Cor}{\textbf{Corollary~}}
\newcommand{\Notation}{\textbf{Notations~}}

\newcommand{\tx}{{\mathfrak t}}
\newcommand{\nx}{{\mathfrak n}}

\newcommand{\Ax}{{\mathfrak A}}
\newcommand{\Bx}{{\mathfrak B}}

\newcommand{\Kc}{{\cal K}}
\newcommand{\Oc}{{\cal O}}
\newcommand{\Xc}{{\cal X}}

\newcommand{\xitl}{{\xi_{\theta, \la}}}
\newcommand{\xitlp}{{\xi_{\theta, \la'}}}
\newcommand{\xitla}{{\xi_{\theta, a\la}}}
\newcommand{\xitpl}{{\xi_{\theta', \la}}}

\newcommand{\chitl}{{\chi_{\theta, \la}}}
\newcommand{\chitlp}{{\chi_{\theta, \la'}}}
\newcommand{\chitla}{{\chi_{\theta, a\la}}}
\newcommand{\chitpl}{{\chi_{\theta', \la}}}

\newcommand{\Jlr}{{J_{\la,\mathrm{right}}}}

\newcommand{\Ch}{{{\mathfrak C}{\mathfrak h}}}
\newcommand{\al}{{\alpha}}
\newcommand{\De}{{\Delta}}
\newcommand{\la}{{\lambda}}
\newcommand{\La}{{\Lambda}}
\newcommand{\tG}{{\widetilde{G}}}

\newcommand{\tGe}{{\widetilde{G_e}}}

\newcommand{\UTn}{{\mathrm{UT}}(n,\Fq)}
\newcommand{\Tn}{{\mathrm{T}}(n,\Fq)}
\newcommand{\Irr}{{\mathrm{Irr}}}

\newcommand{\Fq}{{\Bbb F}_q}
\newcommand{\Cb}{{\Bbb C}}

\newcommand{\row}{{\mathrm{row}}}
\newcommand{\col}{{\mathrm{col}}}

\newcommand{\ind}{{\mathrm{ind}}}
\newcommand{\rank}{{\mathrm{rank}}}

\renewcommand{\leq}{\leqslant}

\title{Supercharacter theory for groups of invertible elements of reduced algebras}
\author{A.N.Panov}
\begin{document}
\Large

\date{}
 \maketitle

\begin{abstract}
 We construct a supercharacter theory for the group of invertible elements of a reduced algebra. For the case of the triangular group, we calculate the  values of   supercharacters on superclasses.\\
\end{abstract}

\section{Introduction}

Traditionally the  problem of classification of irreducible representations is considered to be the main problem in representation theory of finite groups. However, for some groups like the unitriangular $\UTn$, the problem of classification of irreducible representations  appears to be a "wild" problem.
  In 2008 in the paper  \cite{DI}, P.Diaconis and  I.M.Isaacs proposed to replace the problem of classification of irreducible representations by the problem  of construction of a supercharacter theory.  In this theory supercharacters and superclasses play the similar  role as the irreducible characters and classes of conjugate elements in usual theory.
In the same paper    \cite{DI}, the authors presented the supercharacter theory for the algebra groups, i.e. the groups of the type  $1+J$, where  $J$ is a nilpotent  associative finite dimensional  algebra over a finite field. The  main example is the theory of basic  characters for the unitriangular group constructed by C.A.M.Andr\'e ~\cite{A1,A2,A3}  and Ning Yan~ \cite{Yan}.

Notice that a priori there exists several supercharacter theories for a given group.
 The most coarse one has only two supercharacters and two superclasses, and the most precise one ia the theory of irreducible representations.  The  main goal is to construct a supercharacter theory that gives  the most accurate  approximation  of the theory of irreducible characters. The mentioned above theory for the algebra groups is optimal as far as nobody knows the better one.  Since 2008  many papers have been devoted to  the  theory of supercharacters; see the bibliography of the paper  ~\cite{VERY}.  Observe a few papers related to the topic of this article. In the paper  \cite{H}, the general supercharacter theory for semidirect product were presented. However the proposed construction produces rather coarse supercharacter theory, which can be improved in examples.  In the case of semidirect product with abelian normal subgroup,  there are exist more precise  supercharacter theories \cite{SA,Lang}.

In the present paper, we construct  a supercharacter theory for the group of invertible elements of a reduced algebra (theorem \ref{main}).
 This group is semidirect product of the abelian group  $H$ and the algebra group  $1+J$, where  $J$ is the radical of a given algebra  $A$.  The constructed theory is more accurate as in \cite{H}.
 On the other hand, the presented supercharacter theory looks like a continuation of the mentioned above  theory for the algebra groups \cite{DI}. As an example, we obtain a supercharacter theory for  the triangular group  $\Tn$ (theorem \ref{triang}); the values of supercharacters on superclasses are calculated (theorem \ref{vvv}).

 Recall the main definitions of the paper \cite{DI}. Let $G$ be a finite group and $$\Ch = \{\chi_\al|~ \al\in \Ax\}$$ be a system of characters (representations) of the group  $G$. \\
 \Def\Num. The system of characters $\Ch$ defines a supercharacter theory  of the group  $G$ if the characters $\Ch$ are pairwise disjoint and there exists a partition $$ \Kc = \{K_\beta|~ \beta\in\Bx\}$$ of the group  $G$ obeying the conditions:\\
S1) $|\Ax|= |\Bx|$;\\
S2) ~ each character $\chi_\al$ is constant on each subset  $K_\beta$;\\
S3) ~ $\{1\} \in \Kc$ (hereinafter  $1$ is the unit element).
\\
We refer to each  character from $\Ch$ as a {\it supercharacter} and to each subset from  $\Kc$ as a {\it superclass}.

The easily verified condition S3) is very important; one can see this from the following statement from  \cite{DI}. For the readers convenience we present it here with the proof.

Denote by  $\Xc$ the system of subsets  $\{X_\al\}$ in  $\Irr(G)$, where each  $X_\al$ is a subset of irreducible constituents of  the character  $\chi_\al$. \\
\Prop\Num~ \cite[Lemma 2.1]{DI}. Suppose that the system of disjoint characters  $\Ch$ and the system of subsets $\Kc$ obey the conditions S1) and S2).
Then the condition  S3) is equivalent to the following condition \\
S4)~ the system of subsets $\Xc$ defines a partition of  $\Irr(G)$ and each character $\chi_\al$ equals to the character  $$\sigma_\al = \sum_{\psi\in X_\al} \psi(1)\psi $$ up to a constant factor.\\
\Proof.
Suppose that the condition  S3) holds. Every system of disjoint characters  is linear independent. Since $|\Ax|= |\Bx|$, the system $\Ch$ is a basis in the space of all  complex-valued functions over $G$ constant on the subsets from $\Kc$. The regular character  $\rho(g)$ equals to $|G|$, if  $g=1$, and equals to $0$, if $g\ne 1$.
The condition S3) implies that  $\rho(g)$ is constant on the subsets from  $\Kc$. We have
$$\rho = \sum a_\al\chi_\al,~~ \mbox{where}~~ a_\al\in \Cb.$$ On the other hand, each irreducible character is a constituent of  $\rho$ with multiplicity equals to its dimension. Therefore, every irreducible character is a constituent of exactly one  $\chi_\al$ and $\chi_\al = \frac{1}{a_\al}\cdot \sigma_\al$. This proves  S4).

Suppose that the condition S4) holds. The unite of the group belongs to one of the subsets from $\Kc$, say $K_1$. The condition S2) implies  $\chi_\al(g)=\chi_\al(e)$  for any $g\in K_1$ and $\al\in\Ax$. By  S4), we have $\rho(g) = \rho(1)$; this implies  $g=1$ and $K_1=\{1\}$.~$\Box$ \\
\Cor\Num. The system of supercharacters is uniquely determined up to constant factors by the partition  $\Irr(G) = \bigcup X_\al$.

Notice also some important properties of the systems  $\Ch$ and $\Kc$.\\
\Prop\Num ~\cite[Теорема  2.2.]{DI}. \\
1) Each superclass is a union of the classes of conjugate elements.\\
2) The partition  $\Kc$ is uniquely defined by the partition $\Xc$ and vice versa.\\
3) The principal character is a supercharacter (up to a constant factor).

\section{Group of invertible elements of  reduced algebra  }

Let  $A$ be  an unital associative  finite dimensional algebra defined over a finite field   $\Fq$ of $q$ elements. Our main goal is to construct a supercharacter theory for the group  $G=A^*$ of invertible elements of the algebra $A$.

It is known  that, inthe case of a finite field, there exists a complimentary subalgebra $S$ to the radical $J=J(A)$ of the algebra $A$, i.e.  $A = S\oplus J$ (see \cite[\S 11.6]{Pi}). The subalgebra  $S$ is semisimple and is isomorphic to  $A/J$.
The algebra  $A$ is {\it reduced} if the algebra  $A/J$ is commutative (i.e.  $S$ is commutative).
Then  $S$ is a direct sum of fields  \cite[\S 13.6]{Pi}.  Therefore  $$ S=k_1e_1\oplus\ldots\oplus k_ne_n,$$  where $\{e_1,\ldots,e_n\}$  is the system of primitive idempotents, and  $k_1,\ldots, k_n$ are the finite extensions of the field $\Fq$.   Two elements from  $S$ are  {\it associated} if they differ by a factor from  $H$. Any element from  $S$ is associated with some idempotent  $f=e_{i_1}+\ldots+e_{i_m}$  from  $S$.

The group  $G$ is a semidirect product of the subgroup  $H=S^*$ of the invertible elements of the subalgebra  $S$ and the normal subgroup   $N=1+J$. We construct  the group $\tG$ that consists of all triples   $\tau=(t,a,b)$, where $t\in H$ and $a,b\in N$, equipped  with operation  $$((t_1,a_1,b_1)\cdot (t_2,a_2,b_2) = (t_1t_2,~ t_2^{-1}a_1t_2a_2, ~ t_2^{-1}b_1t_2b_2).$$
We define the representation of the group $\tG$ in $J$  by the formula
$$ \rho(\tau)(x) = taxb^{-1}t^{-1}.$$
 One can define  $\tG$-representation in the dual space $J^*$ in the natural way $$\rho^*(\tau)\la(x) = \la(\rho(\tau^{-1})(x)).$$   One can define also left and right $G$-representations in  $J^*$ by formulas  $b\la(x)= \la(xb)$ and $\la a(x) = \la(ax)$. Then $\rho(\tau)(\la) =
tb\la  a^{-1}t^{-1}$.

 For every idempotent  $e\in A$, denote by  $A_e$ the subalgebra $eAe$. The subalgebra  $J_e=eJe$ is the radical in  $A_e$. Denote $e'=1-e$. The Peirce decomposition take place
 $$J= eJe\oplus eJe'\oplus e'Je\oplus e'Je'.$$  We identify  the dual space $J_e^*$ with the subspace in $J^*$ of all linear forms equal to zero on all components of the Peirce decomposition apart from the first one. \\
\Def\Num. \\
1) We refer to    $x\in J$
 as {\it singular} element  if the exists an idempotent  $e\in A$, ~ $e\ne 1$, such that $x\in J_e$. Otherwise we refer to   $x$ as {\it regular}.  \\
2) We refer to $\la\in J^*$
as {\it singular} element  if there exists an idempotent  $e\in A$, ~ $e\ne 1$,  such that  $\la\in J_e^*$. Otherwise we refer to   $\la$ as  {\it regular}.\\
\Lemma\Num\label{ccc}. The following conditions are equivalent:
1) an element  $x\in J$ (respectively,  $\la\in J^*)$ is  singular;
2) there exists  $c\in A\setminus J$ such that  $cx = xc = 0$ (respectively, $c\la = \la c= 0$).
\\
\Proof. We shall prove for  $x\in J$. The case  $\la\in J^*$ is considered similarly.

Notice that the condition  $x\in J_e$ is equivalent to  $x\in J$ and $e'x=xe'=0$.
So the element  $x$ is  singular if there exists an idempotent   $f\ne 0$ such that $fx=xf=0$.

It is evident that   1)  implies  2). Let us prove that  2)  implies  1). Let  $c$ belongs to  $A\setminus J$ and  $cx = xc = 0$.
Suppose that  $x\ne 0$ (the zero element is  singular).
The element $c$ generates the  commutative subalgebra $C=\Fq[c]$ that has the ideal  $I=c\Fq[c]$. Since  $x\ne 0$, the codimension of $I$ in $C$ equals to 1.

Suppose that  $C$ is a local algebra. In this case,  $C$ has a unique ideal of codimension 1 and it consists of nilpotents; in particular, $c$ is nilpotent. On the other hand,   as $c$ does not belong to $ J$, it is  not nilpotent. Therefore, $C$ is not a local algebra.

The algebra  $C$ contains idempotent  $f$ doesn't equal to 0 or 1.
 If $f\in I$, then the statement is proved.
 Let $f\notin I$. There exists  $\al\in \Fq$ such that $\al-f\in I$. Then $$(\al-f)^2=\al^2-2\al f+f^2=\al^2+(1-2\al)f\in I.$$
Easy to see that \begin{equation}
\left|\begin{array}{cc} \al&-1\\
\al^2&1-2\al\end{array}\right| = \al(1-\al).
\end{equation}
 If $\al\ne 0,1$, then $1\in I$ and  $I=C$; this contradicts to  $\mathrm{codim}(I, C)=1$.

 If $\al=0$, then $f\in I$; this contradicts to choice of  $f$. Hence  $\al=1$. In this case, the idempotent  $1-f\ne 0$ belongs to  $I$. This proves the statement. $\Box$\\
\Prop\Num. If the element  $x\in J$ (respectively, $\la\in J^*$) is  singular,  then every element of its   $\rho(\tG)$-orbit   (respectively, $\rho^*(\tG)$-orbit) is also singular. Similarly for  regular elements.\\
\Proof. We shall prove for  $x\in J$. The case  $\la\in J^*$ is treated similarly. Let $x$ be a  singular element. There exists an idempotent $f\in A$, ~ $f\ne 0$, such that  $fx=xf=0$. It is sufficient to prove that $ax$ (respectively, $xa$) is  singular  for any $a\in N$.  Then  for  $c=fa^{-1}$, we have $c\notin J$ and $c(ax)=(ax)c=0$. The Lemma \ref{ccc}
implies that the element  $ax$ is singular. $\Box$\\
\Cor\Num\label{einS}. For any  singular  $\tG$-orbit, $\Oc\subset J$ (respectively, $\Oc^*\subset J^*$) there exists an idempotent  $e\in S$
   such that $\Oc\cap J_e\ne \varnothing$ (respectively, $\Oc^*\cap J_e^*\ne \varnothing$).\\
\Proof. Well known that every idempotent is conjugate to an idempotent from  $S$ (see for example \cite[теорема 4.1]{DK}). $\Box$

Denote by $\tGe$ the group of invertible elements of $A_e$. Identifying  $\tilde{g_e}\in \tGe$ with   $ \tilde{g_e}+e'\in \tG$, where $e'=1-e$,
one can treat  $\tGe$ as a subgroup of  $\tG$.\\
\Lemma\Num\label{Oc}. Let  $e,f$ be two idempotents from  $S$.  Then\\
1)~ if an  $\tG$-orbit $\Oc$ has nonempty intersection with both  $J_e$ and $J_f$, then  $\Oc$ has  nonempty intersection with  $J_{ef}$; \\
2)  two elements   $x $ and $y$ from $J_e$ belong to a common  $\tG$-orbit  if and only if $x$ and $y$ belong to a common  $\tGe$-orbit. \\
\Proof.  \\
1) Let  $x\in \Oc\cap J_e$, ~$y\in \Oc\cap J_f $ (i.e.  $ex=xe=x$, ~
 $fy=yf=y$, and $y=haxbh^{-1}$, where $h\in H$ and $a,b\in N$). Consider the element  $z=eye$. Since  $ef=fe$, we have  $(ef)z=z(ef)=z$ and $z=eye= ehaxbh^{-1}e = heaxbeh^{-1}$.
The element   $a_1= ea+e'$, where  $e'=1-e$, obeys  $a_1x=eax+(1-e)x=eax$. As $a_1 = 1\bmod J$, we obtain $a_1\in N$.
Similarly,  $b_1=be+ e'$ belongs to  $N$ and $xb_1=xbe+x(1-e) = xbe$. Then $z=ha_1xb_1h^{-1}$, where $a_1,~b_1\in N$; that is
$z$ lies in the common orbit  $\Oc$.\\
2) Let  $x,y\in J_e$ (i.e. $ex=xe=x$ and $ey=ye=y$). It is evident that if $x,y$ belong to a common  $\tGe$-orbit, then  $x,y$ belong to a common  $\tG$-orbit. Let us prove the inverse statement.

  Let $x,y$ belong to the common $\tG$-orbit, i.e.  $y=haxbh^{-1}$, where $a,b\in N$ and $h\in H$. Then
$$y=eye=ehaxbh^{-1}e = (he)a(exe)beh^{-1} = (he)(eae)x(ebe)(h^{-1}e).$$ The element  $he=eh$ is invertible in $A_e$ and its inverse coincides with $h^{-1}e$. The elements  $eae$ and $ebe$ belong to  $N_e$. This proves that  $x$ and $y$ lie in a common  $\tGe$-orbit. $\Box$\\
\Cor\Num\label{eee}. For every   $x\in J$ there exists a unique idempotent $e\in S$ such that  $\Oc(x)\cap J_e\ne \varnothing$ and   $\Oc(x)\cap J_e$ is a regular  $\tGe$-orbit in  $J_e$. \\
\Proof. The proof follows from statements 1) and  2) of the previous  lemma.
$\Box$
\\
The following properties of orbits in  $J^*$ are proved similarly.\\
\Lemma\Num\label{Oc*}. Let  $e,f$ be idempotents in  $S$. Then\\
1)~ if an  $\tG$-orbit $\Oc^*$ has nonempty intersection with both  $J_e^*$ and $J_f^*$, then  $\Oc^*$ has  nonempty intersection with  $J^*_{ef}$; \\
2) two elements   $\la $ and $\mu$ from $J_e^*$ belong to a common  $\tG$-orbit  if and only if $\la$ and $\mu$ belong to a common  $\tGe$-orbit.\\
\Cor\Num\label{eeeconj}. For any  singular  $\la\in J^*$ there exists a unique idempotent  $e\in S$ such that $\Oc(\la)\cap J^*_e\ne \varnothing$ and   $\Oc(\la)\cap J^*_e$ is a  regular  $\tGe$-orbit in  $J^*_e$.\\
\Notation\Num: \\
$\Oc(J)$ is the set  of all  $\tG$-orbits in $J$;\\
$n(J)$ is the number of all  $\tG$-orbits in  $ J$;\\
$\Oc(J_e)$ is the set of all  $\tG$-orbits of the type $\Oc(x)$, where  $x\in J_e$; \\
$n(J_e) = |\Oc(J_e)|$ ( by Lemma  \ref{Oc}  the number   $n(J_e)$ coincides with the number  of all  $\tGe$-orbits in $J_e$);\\
$n_E(J)$ is the number of all  regular $\tG$-orbits in  $ J$;\\
$\Oc(J^*)$,~ $n(J^*)$,~ $\Oc(J^*_e)$,~ $n(J^*_e)$,~ $n_E(J^*)$ are defined similarly.
\\
 \Cor\Num\label{cap}.~ 1) $\Oc(J_e)\cap \Oc(J_f) = \Oc(J_{ef})$,~~ 2) ~$\Oc(J_e^*)\cap \Oc(J_f^*) = \Oc(J_{ef}^*)$.
\Prop\Num\label{ereg}. The number of all regular  $\tG$-orbits  in  $J$ coincides with the number  of all  regular  $\tG$-orbits  in  $J^*$ (i.e.  $n_E(J) = n_E(J^*)$).\\
\Proof.
  Let $\{e_1,\ldots, e_n\}$ be a system of all primitive idempotents of the algebra  $S$. Denote $e'_i = 1-e_i$.
For each idempotent  $f\in S$, we denote by $l(f)$ the number of factors in the decomposition  $f=e'_{i_1}\cdots e'_{i_l}$. For  $f=1$, we take  $l(f)=0$.
The Corollary   \ref{cap} implies that
$\Oc(J_{f})$  coincides with the intersection  $\bigcap \Oc(J_{\phi})$, where $\phi $ runs through the set of idempotents   $\{e'_{i_1},\ldots, e'_{i_l}\}$.

The set of all singular orbits  in   $J$ is a union of all $\Oc(J_{e'_i})$, where $1\leq i\leq n$. Hence
$$ n_E(J) = \bigoplus (-1)^{l(f)} n(J_f),$$
where  $f$ runs through the set of all idempotents of  $S$. The similar formula is true for the orbits in  $J^*$:
$$ n_E(J^*) = \bigoplus (-1)^{l(f)} n(J^*_f), $$
where $f$ is also runs through the set of all idempotents of $S$.
For any finite subgroup of linear operators in a linear space $V$ over a finite field, the number of orbits in  $V$ coincides with the number of orbits in the dual space  $V^*$ (см.\cite[Lemma 4.1]{DI}). For the orbits of the group  $\tilde{G_f}$ in  $J_f$  and $J^*_f$, we conclude that $n(J_f) = n(J^*_f)$. Therefore  $n_E(J)=n_E(J^*)$. $\Box$

\section{Superclasses}

Let us define the action of the group  $\tG$ on  $G$
as follows
\begin{equation}\label{RRR}
    R_\tau(g) = 1+ ta(g-1)b^{-1}t^{-1}, \quad \mbox{где}\quad  \tau = (t, a, b), ~ t\in H,~ a,b\in N.
\end{equation}
Notice that for  $g=1+x\in N$, we obtain $R_\tau(1+x) = 1+\rho_\tau(x)$.

A {\it superclass} $\Kc(g)$ is an  $\tG$-orbit of  the element $g\in G$. Every element $g\in G$ can be presented in the form  $g=h+x$, where  $h\in H$ and $x\in J$. Easy to see that if  $g=h+x$ and $g'=h'+x'$ lie in a common superclass, then $h=h'$.
\\
\Theorem\Num\label{superclass}. Let   $g=h+x$, where $h\in H$ and $x\in J$. Let  $\Kc=\Kc(g)$. Let  $f$ be the idempotent associated with  $h-1$. Then \\
1) there exists an idempotent  $e\in S$, ~ $e\perp f$, and a regular element  $y\in J_e$ such that  $h+y$ belongs to the superclass  $\Kc$; denote by  $\omega(y)$ the orbit of the element  $y$ in  $J_e$ with respect to the group $\tGe$; \\
2) the quadruple   $\beta = (e,f,h, \omega)$ uniquely determines the superclass $\Kc$. \\
\Proof.\\
{\bf Item 1}. Let us prove that there exists an element  $h+y\in \Kc$ such that  $y\in J_{f'}$, where $f'=1-f$.

 Denote $s=h-1$. It is sufficient to prove that there exist $a,b\in N$ such that
$s+y=a(s+x)b$, where $yf = fy = 0$.
Since  $s$ is associated with the idempotent  $f$, it is sufficient to prove that for  any $x\in J$ there exist
$a,b\in N$  and $y\in J$ such that  $f+y=a(f+x)b$ and $yf=fy=0$.
\\
1) Let us show that there exists  $u\in J$ such that $f+u=a(f+x)$ and $uf=0$.
 Take $a=(1+x)^{-1}\in N$. Then
 $$
 uf=\left((1+x)^{-1}(f+x)-f\right)f= ((1+x)^{-1}(f+xf)-f = 0.$$
 2) Let  $u$ be as in 1). Let us show that there exists $y\in J$ such that  $f+y = (f+u)b$,~ $ b\in N$, and $yf=fy=0$.  Take   $b=(1-fu)\in N$. Then
          $$y=(f+u)(1-fu)-f=(1-f)u.$$ Hence $fy=yf=0$.\\
 {\bf Item 2}.
In this item, we shall prove that the elements $h+y$ and $h+y'$, where $y, y'\in J_{f'}$, belong to a common superclass if and only if  $y$ and $y'$  belong to a common $\tG_{f'}$-orbit in   $J_{f'}$.

It is obvious that if elements $y$ and  $y'$ belong to a  common $\tG_{f'}$-orbit, then $h+y$ and  $h+y'$ belong to a common superclass.  Let us prove the inverse statement.

Let  $h+y$ and $h+y'$ belong to a common supercalss. Then there exist  $t\in H$, ~ $u,v\in J$ such that
\begin{equation}\label{hhh}
 h-1 + y' = t(1+u)(h-1 + y)(1+v)^{-1}t^{-1}
\end{equation}
Then
\begin{equation}\label{huv}
 (h-1 + y')t(1+v) = t(1+u)(h-1 + y)
\end{equation}
Multiply (\ref{huv}) to the left and to the right by $f'$. Applying  $f'(h-1)=(h-1)f'=0$, ~ $f'y=yf'=y$, ~  $f'y'=y'f'=y'$ and $f't=tf'$, we obtain
$$ y't(1+f'vf') =  t(1+f'uf')y.$$
That is  $y$ and $y'$ belong to a common  $\tG_{f'}$-orbit in  $J_{f'}$.\\
{\bf Item  3}. For any  $\tG_{f'}$-orbit $\omega$ in  $J_{f'}$, there exists a unique idempotent  $e<f'$ such that  $\omega\cap J_e$ is a  regular  $\tGe$-orbit in $J_e$ (see Corollary  \ref{eee}).
One can consider that $y$ is a  regular element in $J_e$. This proves the statement  1). The Lemma \ref{Oc} implies the statement 2) of the theorem. $\Box$

Denote by  $\Bx$ the set of all quadruples  $\beta = (e,f,h, \omega)$, where  $e, f$ is a pair of edempotents  $S$, where $e\perp f$, ~ the element   $h\in H$ and  $h-1$ is associated with  $f$, and $\omega$ is a regular  $\tGe$-orbit in  $J_e$.\\
\Notation\Num\label{nota}.\\
For any  $\beta\in\Bx$, we denote by  $\Kc_\beta$ the superclass of the element  $g=h+x$, where $x\in \omega$;\\
$S_i = e_iS = Se_i; \quad H_i=\{1-e_i+s_i, ~~ \mbox{where}~~ s_i\in S_i^*\};$\\
\\
$ H_f =\prod_{e_i\leq f} H_i;\quad m(f) = \sum_{e_i\leq f} (|H_i|-1) $.\\
\\
\Cor\Num\label{supcl}. The superclasses  $\{\Kc_\beta|~ \beta\in \Bx\}$ form a partition of the group  $G$. The number of superclasses equals to
\begin{equation}\label{BxBx}
|\Bx| = \sum_{e\perp f}(n_E(J_e)+m(f)).
\end{equation}

\section{Supercharacters}

Let $e,f$ be idempotents from  $S$ and $e\perp f$. Let $\la$ be a regular element in  в $J_e^*$. As above  $e'=1-e$ and $$ H_{e'} = \prod_{e_i\leq e'} H_i.$$
Notations:
$$J_{\la,\mathrm{right}} = \{x\in J|~ \la x =0\},\quad\quad  N_{\la,\mathrm{right}} = \{ a\in N|~ \la a =\la\}.$$
Notice that  $N_{\la,\mathrm{right}} = 1+ J_{\la,\mathrm{right}}$.
We define similarly $H_{\la,\mathrm{right}}$ and $H_{\la,\mathrm{left}}$.

Easy to see that  $H_{e'}$ is contained in  $H_{\la,\mathrm{right}}$ and $H_{\la,\mathrm{left}}$. Because   $\la$ is  a regular element in  $J_e^*$, we obtain the equality
\begin{equation}\label{}
    H_{e'} = H_{\la,\mathrm{right}}\cap H_{\la,\mathrm{left}}.
\end{equation}
Consider the subgroup   $G_\la = H_{e'}\cdot N_{\la,\mathrm{right}}$. The subgroup  $G_\la$ is a semidirect product of  $H_{e'}$ and the normal subgroup  $N_{\la,\mathrm{right}}$.  Every element  $g\in G_\la$ can be presented in the form  $g=h+x$, where $h\in H_{e'}$ and $x\in J_{\la,\mathrm{right}}$.

Fix a nontrivial character  $c\to  \varepsilon^{c}$  of the additive group of the field  $\Fq$ with values in the multiplicative group  $\Cb^*$. Let $\theta$ be a linear character
 (one-dimensional representation)   of the subgroup $H_{e'}$.
Define a linear character of the subgroup  $G_\la$ as follows
\begin{equation}\label{}
    \xi_{\theta,\la}(g) = \theta(h)\varepsilon^{\la(x)},
\end{equation}
where $g=h+x$,~ $h\in H_{e'}$ and $x\in J_{\la,\mathrm{right}}$.
Let us show that  $\xi=\xi_{\theta,\la}$ is really a linear character:
$$\xi(gg') = \xi((h+x)(h'+x')) = \xi(hh' + h'x + x'h  + xx') =$$
$$\theta(hh')\varepsilon^{\la(h'x)} \varepsilon^{\la(x'h)}\varepsilon^{\la(xx')} = \theta(h)\theta(h')\varepsilon^{\la(x)} \varepsilon^{\la(x')} = \xi(g)\xi(g').$$
We refer to the  induced character
\begin{equation}\label{indchi}
\chi_{\theta,\la} = \ind(\xi_{\theta,\la}, G_\la, G).
\end{equation}
as a {\it supercharacter}.\\
\Prop\Num\label{supersuper}. The supercharacter   $\chi_{\theta,\la}$ is constant on the superclasses.\\
\Proof. \\
{\bf Item  1}. Let  $g\in G_\la$, ~$a\in N$. Let us show that $g'=1+(g-1)a\in G_\la$ and $\chitl(g')=\chitl(g)$.

The element  $g\in G_\la$ is presented in the form  $g=h+x$, where $h\in H_{e'}$ and $x\in \Jlr$. Let $a=1+u$, where $u\in J$.
Then $g'= h+y$, where $h\in H_{e'}$ and $y=(h-1)u + x+ xu\in J$.  Since $$\la h =\la,\quad \la(h-1)u=0, \quad \la x=0~~\mbox{и}~~ \la xu=0,$$  we have
$\la y = \la (h-1)u +\la x + \la xu = 0$. Hence $y\in J_{\la,\mathrm{right}}$ and $g'\in G_\la$.

We obtain
\begin{equation}\label{xixi}
    \xitl(g') = \theta(h)\varepsilon^{\la((h-1)u+x+xu)}=\theta(h)\varepsilon^{\la(x)}=\xitl(g).
\end{equation}
Denote $\La(g) = \{s\in G|~ s^{-1}gs\in G_\la\}$.

 Let us show that $\La(g') = \La(g)$. Really,  $s^{-1}g's = 1+ (s^{-1}gs-1)s^{-1}as$.
   As  above,  $ s^{-1}gs\in G_\la$ implies $s^{-1}g's\in G_\la$. This proves that  $\La(g)\subset\La(g')$. The reverse inclusion similarly.

The formula  (\ref{xixi}) implies
\begin{equation}
\xitl(s^{-1}g's) = \xitl(s^{-1}g s),
\end{equation}
for  $s\in \La(g)$.
Therefore
\begin{equation}\label{chichi}
\chitl(g') = \frac{1}{|G_\la|}\sum_{s\in \La(g')} \xitl(s^{-1}g's)= \frac{1}{|G_\la|}\sum_{s\in \La(g)} \xitl(s^{-1}gs) = \chitl(g).
\end{equation}
{\bf Item 2}. The characters are constant on the classes of conjugate elements. Hence $\chitl(g_0gg_0^{-1})=\chitl(g)$ for every  $g_0\in G$. From formula (\ref{RRR}) and Item 1, we conclude that
$\chitl(R_\tau(g))=\chitl(g)$ for every  $\tau\in\tG$ and $g\in G_\la$.

If an $\tG$-orbit has empty intersection with  $G_\la$, then is $\chitl$ in zero on it.$\Box$\\
\Prop\Num\label{lalap}.
 Let $\la,~\la'$ be regular elements from $J_e^*$. If $\la, \la'$  belong to a common  $\tGe$-orbit, then $\chitl=\chitlp$.\\
\Proof.  The statement follows from the following items 1 and 2.\\
{\bf Item 1}. Consider the case  $\la'=a\la$, where $a=1+u\in N_e$. Since  $ J_{a\la,\mathrm{right}}
=J_{\la,\mathrm{right}}$, we have $G_{a\la}=G_\la$. For any   $g=h+x\in G_\la$ we obtain
$$ \xitla(g) =   \theta(h)\varepsilon^{a\la(x)} = \theta(h)\varepsilon^{\la(x+xu)}= \theta(h)\varepsilon^{\la(x)}=\xitl(g),$$
where $g=h+x$, ~ $h\in H_{e'}$ and $x\in J_{\la,\mathrm{right}}$.
Hence $\chitl=\chitla$.\\
{\bf Item 2}. Let $\la'= g_0\la g_0^{-1}$, where $g_0\in G_e$. Then $G_{\la'}=g_0G_\la g_0^{-1}$ and
$$\xitlp(g) = \theta(h)\varepsilon^{g_0\la g_0^{-1}(x)} = \theta(h)\varepsilon^{\la(g_0^{-1}xg_0)} = \xitl(g_0^{-1}gg_0).$$
The  corresponding induced representations are equivalent and $\chitl=\chitlp$. $\Box$\\
\Prop\Num\label{thethep}.
 If $\theta\ne\theta'$, then the characters  $\chitl$ and $\chitpl$ are disjoint.\\
\Proof. Denote  $\xi=\xitl$ and $\xi'=\xitpl$. It follows from the Intertwining Number Theorem \cite[теорема 44.5]{CR} that
the characters  $\chitl$ and $\chitpl$ are disjoint if and only if  for every $s\in G$ there exists  $g\in G_\la$ such that  $sgs^{-1}\in G_\la$ and $\xi(sgs^{-1})\ne \xi'(g)$.

It is sufficient to prove that for every $t\in H$ and $a\in N$ there exists  $g\in G_\la$ such that  $aga^{-1}$ and  $tgt^{-1}$ belong to $G_\la$, ~  and
\begin{equation}\label{at}
    \xi(aga^{-1})\ne \xi'(tg t^{-1}).
\end{equation}
{\bf Item 1}.  Let  $a\in N$ and $t\in H$. In this item, we reduce the proof of  (\ref{at}) to the case  $a=1+ev$, where   $v\in J$.

Let us show that there exists  $w\in J$ such that  $a(1+e'w)=1+ev$, where  $v\in J$,~ $e'=1-e$.
Really, let  $a=1+u$,~ $u\in J$. Then  for every  $w\in J$, we have:
$$ a(1+e'w) =  (1+eu+e'u)(1+e'w) = 1+ eu(1+e'w) + e'u + (1+e'u)e'w.$$
It is sufficient to take  $e'w = -(1+e'u)^{-1}e'u$.

  Since $\la e'=0$, we  see  $e'w\in J_{\la,\mathrm{right}}$ and $1+e'w\in N_{\la,\mathrm{right}}$. Moreover,
$$t(1+e'w)t^{-1}= 1+e' twt^{-1}\in N_{\la,\mathrm{right}}\subset G_\la.$$
 We conclude that it is sufficient to prove (\ref{at}) for $a=1+ev$, where  $v\in J$.\\
{\bf Item  2}. Let $t\in H$ and  $a=1+ev=1+eve+eve'$, where $v\in J$.
Let us prove  (\ref{at}) for the case  $eve'\in J_{\la,\mathrm{right}}$.

Notice that  $(1-eve')a=(1-eve')(1+eve+eve') = 1+eve$. Since
$1-eve'\in N_{\la,\mathrm{right}}\subset G_\la$, it is sufficient to prove  (\ref{at})  for $a=1+eve$, where  $v\in J$.

Let $g=h_0$ be an arbitrary element of the subgroup  $H_{e'}\subset G_\la$.
Then  $h_0e=eh_0=e$ and
$$ah_0a^{-1}= (1+eve)h_0(1+eve)^{-1} = h_0(1+eve)(1+eve)^{-1} =h_0\in G_\la.$$
Since $t$ and $h_0$ commutes,  $th_0t^{-1} = h_0\in G_\la$.

By assumption, $\theta\ne\theta'$; there exists  $h_0\in H_{e'}$ such that  $\theta(h_0)\ne\theta'(h_0)$. Taking   $g=h_0$ in (\ref{at}), we obtain
 $$\xi(ah_0a^{-1}) = \xi(h_0) \ne \xi'(h_0) = \xi'(th_0t^{-1}).$$
 {\bf Item 3}. Let  $a=1+x$, where $x=eve+eve'\in J$ and  $eve'\notin J_{\la,\mathrm{right}}$.
 Then $$\la(xe'J) = \la((eve+eve')e'J)=\la((eve')e'J)=\la(eve'(eJ+e'J))=\la(eve'J)\ne 0.$$
Consider the chain of right ideals  $e'J\supset e'J^2\supset\ldots\supset e'J^k=\{0\}.$
There exists a number  $i$ such that  $\la(xe'J^i)\ne 0$ and $\la(xe'J^{i+1})= 0$.
Choose  $y\in e'J^i$ such that  $\la(xy)\ne 0$ and $\la(xyJ) = 0$. The last equality means that   $xy\in J_{\la,\mathrm{right}}$.
Since $\la e'=0$, we conclude  that $y\in J_{\la,\mathrm{right}}$.

Take  $g=1+cy\in N_{\la,\mathrm{right}}$, where $c\in \Fq$. Then
$$ aga^{-1} =(1+x)(1+cy)(1+x)^{-1} =1+(1+x)cy(1+x)^{-1} \in 1+cy + cxy + yJ + xyJ$$
Therefore, $aga^{-1}\in N_{\la,\mathrm{right}}\subset G_\la.$
As $\la(e'J)=0$, we have $\la(y)=0$.
Since  $y,~xy\in J_{\la,\mathrm{right}}$, we have $\la(yJ) = \la(xyJ)=0$.
We obtain
$\xi(aga^{-1}) = \varepsilon^{c\la(xy)}$.

Since
$y\in e'J$, we see  $tyt^{-1}\in e'J$ and  $$tgt^{-1} = 1+tyt^{-1}\in N_{\la,\mathrm{right}}\subset G_\la.$$
We  have  $\xi'(tgt^{-1}) =\varepsilon^{c\la(tyt^{-1})}=1$.

 As
 $\la(xy)\ne 0$, there exists  $c\in \Fq$ such that $\varepsilon^{c\la(xy)}\ne 1$. This proves  (\ref{at}). $\Box$

Let  $e\perp f$ be a pair of orthogonal idempotents. We say that a linear  character $\theta$ of the subgroup  $H_{e'}$ is {\it associated } with idempotent $f$ if for any primitive idempotent  $e_i\leq e'$, the following conditions hold:\\
1) if $e_i\leq f$, then  $\mathrm{Res}_{H_i}\theta \ne 1$,\\
2) if $e_i\nleqslant f$, then $\mathrm{Res}_{H_i}\theta = 1$.\\
The number of characters associated with $f$ equals to $m(f)$ (see Notations \ref{nota}).

 Denote by $\Ax$ the set of all quadruples  $\al = (e,f,\theta, \omega^*)$, where $e, f$ is a pair of idempotents from $S$, ~  $e\perp f$, ~ $\theta$ is a linear character of the subgroup  $H_{e'}$, associated with $f$, and       $\omega^*$ is a  regular  $\tGe$-orbit in  $J^*_e$.
 For any $\al\in\Ax$, we denote by  $\chi_\al$ the supercharacter equals to  $\chitl$, where  $ \la\in\omega^*$ (see Proposition \ref{lalap}).\\
\Prop\Num\label{dva}. The characters $\{\chi_\al|~ \al\in \Ax\}$  are pairwise disjoint.\\
\Proof.   First recall the construction of supercharacters for the algebra group  $N=1+J$ (see \cite{DI}). For any $\mu\in J^*$, the supercharacter  $\chi_\mu$ of the group $N$ is defined as an induced character from the character  $$\xi_\mu(1+x)=\varepsilon^{\mu(x)},$$
 of the subgroup $N_{\la,\mathrm{right}}$. Two supercharacters $\chi_\mu$ and $\chi_{\mu'}$ are equal if and only if  $\mu$  and  $\mu'$ belong to a common $N\times N$-orbit. Otherwise, the supercharacters  $\chi_\mu$  and $\chi_{\mu'}$ are disjoint ~\cite[Theorem 5.5]{DI}.

Applying the formula  of decomposition of a restriction of
induced representation on  a subgroup (see \cite[теорема 44.2]{CR}, \cite[proposition 22]{Serr}), we have \begin{equation}\label{res}
   \mathrm{Res}_N (\chitl) = \sum_{t\in H_e} \chi_{Ad^*_t\la}.
   \end{equation}

 Turn to the proof of the statement. Let $\al\ne \bar{\al}$, where $\al = (e,f,\theta, \omega^*)$ and  $\al = (\bar{e},\bar{f},\bar{\theta}, \bar{\omega}^*)$.
 If $e\ne \bar{e}$, then $\la $ and $\bar{\la}$ belong to different $\tG$-orbits (see
 Corollary \ref{eeeconj}).
By formula  (\ref{res}),
the restrictions of  $\chi_\al$ and $\chi_{\bar{\al}}$  on $N$ are disjoint.
It follows that $\chi_\al$ and $\chi_{\bar{\al}}$  are disjoint.

Similarly, we treat the case $e=\bar{e}$, ~$f=\bar{f}$,~ $\theta=\bar{\theta}$ and
$\omega^*\ne \bar{\omega}^*$ (see Lemma \ref{Oc*}).

It remains to consider the case  $e=\bar{e}$, ~ $\omega^*= \bar{\omega}^*$ ,
$f\ne\bar{f}$ or $\theta\ne\bar{\theta}$. The proof follows from Proposition \ref{thethep}. $\Box$\\
  \Theorem\Num\label{main}. The systems of supercharacters $ \{\chi_\al|~ \al\in \Ax\}$ and  superclasses  $\{\Kc_\beta|~ \beta\in \Bx\}$  define a supercharacter theory on the group  $G$.\\
 \Proof. The supercharacters are disjoint according to the previous proposition. The superclasses form a partition of the group  $G$ (see Corollary  \ref{supcl}). The number of supercharacters  is equal to
\begin{equation}\label{AxAx}
|\Ax| = \sum_{e\perp f}(n_E(J^*_e)+m(f)).
\end{equation}
By  formula (\ref{BxBx}) and Proposition  \ref{ereg}, we conclude that  $|\Ax|=|\Bx|$; this proves  S1). The supercharacters are constant on superclasses  (see Proposition \ref{supersuper}); this proves  S2). Finally,  $\{1\}$ is the superclass $K(g)$ for $g=1$; this proves  S3). $\Box$

\section{Supercharacter theory for the triangular group}

 Consider the algebra  $A=\tx(n,\Fq)$ of all  $n\times n$-matrices with entries from the field $\Fq$ and with zeros
 under the diagonal. The {\it triangular group}  $G=\Tn$ is the group of invertible elements of the algebra  $A$. The radical  $J$ coincides with the subalgebra  $\nx(n,\Fq)$ of all triangular matrices with zeros on the diagonal. The complementary subalgebra  $S$ is the subalgebra of diagonal matrices; ~$H=S^*$ is the subgroup of diagonal matrices. In this section, we specialize the constructed supercharacter theory for the triangular group.

 Simplifying the language, we refer to a pair
  $(i,j)$, where $i,j$ are  positive integers  $1\leq i<j\leq n$, as a  {\it
positive root}. Denote by  $\De$ the set of all positive roots. We  refer to $i=\row(\gamma)$ (respectively,  $j=\col(\gamma)$), as a number of row (respectively, column) of the root $\gamma=(i,j)$.
 Following C.A.M.Andr\'e, we introduce definition of a basic subset. A  {\it basic subset}  is a subset   $D\subset\De$
 that has no more than one positive root in any row and any column ~ \cite{A1,A2}.

Let $\{E_{ij}|~1\leq i<j\leq n \}$ be a basis of matrix units of the radical  $J$. To any basic subset $D$, we attach the element  $$x_D=\sum_{\gamma\in D}  E_{ij},$$ in the radical  $J$.
 Denote by  $\Oc_D$  the  $\tG$-orbit of the element $x_D$. The orbits $\{\Oc_D\}$ form a partition of $J$ (see \cite{A2}); therefore,  $\{1+\Oc_D\}$  form a partition of $N$.

Let  $\{E^*_{ij}\}$ be the dual basis with respect to the basis $\{E_{ij}\}$. Next, we attach to $D$ the linear form
$$\la_D=\sum_{\gamma\in D}  E^*_{ij}$$
and its $\tG$-orbit $\Oc^*_D $.
The orbits  $\{\Oc^*_D\}$ form a partition of  $J^*$~\cite{A1}.\\
\Lemma\Num\label{nnn}. The orbit  $\Oc_D$ (respectively,  $\Oc^*_D$) is   regular if and only if  $\row(D)\cup\col(D)=[1,n]$.\\
\Proof. We shall prove that the  orbit  $\Oc_D$  is   singular if and only if  $\row(D)\cup\col(D)\ne [1,n]$.

Let $\row(D)\cup\col(D)\ne[1,n]$, then the element $x_D$ belongs to  $J_e$, where
$e=\sum E_{ii}$ and  $i$ runs through $ \row(D)\cup\col(D)$.
   This proves that $x_D$ (and the orbit $\Oc_D$) is  singular.

On the other hand, suppose that  $\Oc_D$  is a singular orbit; then there exists an idempotent  $e\in S$, ~ $e\ne 1$, such that  $\Oc_D\cap J_e\ne \varnothing$ (see Corollary \ref{einS}). The algebra  $J_e$ is isomorphic to $\nx(k,\Fq)$ for some   $k<n$. The intersection $\Oc_D\cap J_e$ is a  $\tGe$-orbit (see Lemma ~\ref{Oc}); there exists $x_{D'}\in \Oc_D\cap J_e$ for some  basic subset $D'$. Since $x_D$ and $x_{D'}$ belong to a common  $\tG$-orbit, we conclude $D=D'$. Therefore
$$\row(D)\cup\col(D)=\row(D')\cup\col(D')\ne [1,n].~\Box$$

 Let $h=\mathrm{diag}(h_1,\ldots,h_n)$ be an  element of  $H$, obeying  $h_i=1$ for each $i\in \row(D)\cup\col(D)$. Let  $g_{h,D} = h+x_D\in G$ and $K_{h,D}$ be a $\tG$-orbit of the element  $g_{h,D}$ in the group  $G$.

Let $\theta$ be a linear character  of the group  $H$, obeying
$\mathrm{Res}_{H_i}(\theta) = 1$ for each $i\in \row(D)\cup\col(D)$.
For  $\la_D$ and  $\theta$ similarly to   (\ref{indchi}), we define the induced character  $$\chi_{\theta,D} = \chi_{\theta,\la_D} .$$
\Theorem\Num\label{triang}. The systems of characters  $\{\chi_{\theta,D}\}$ and subsets  $\{K_{h,D}\}$ define a supercharacter theory on the group  $G=\Tn$.\\
\Proof. The proof follows from Lemma \ref{nnn} and Theorem  \ref{main}. $\Box$

Calculate the  values of supercharacters  $\chi_{\theta,D}$ on superclases $K_{h,D'}$. We need new notations.

For each positive root  $\gamma=(i,j)$, ~$1\leq i<j\leq n$, we denote
$$ \Delta'(\gamma) =\{(i,k)| ~i<k<j\},\quad\quad\Delta''(\gamma) = \{(k,j)|~
i< k<j\}.$$
 The number of elements of both subsets coincides and equals to  $j-i-1$.

 Take   $ \delta' (D,D') = 0$ if there exists  $\gamma\in D$ and $\gamma\in D'$ such that  $\gamma'\in\De'(\gamma)$. Otherwise we take $\delta'(D,D')=1$.
Similarly we define  $ \delta''(D,D')$.

 We take $\delta_0(D,h) = 1$, if  $h_i=1$ for each  $i\in \row(D)\cup\col(D)$;
otherwise   $\delta_0(D,h) = 0$. Denote
\begin{equation}\label{delta}
    \delta(D,h,D') =  \delta'(D,D') \delta''(D,D')\delta_0(D,h).
\end{equation}

For each positive root  $\gamma=(i,j)$, ~$1\leq i<j\leq n$, we denote by
$P(\gamma)$ the submatrix of the matrix  $g-1$ with system of rows and columns $[i+1,j-1]$. Denote by    $m(\gamma, h, D')=\mathrm{corank}P(\gamma)$.  Since the matrix $P_\gamma$ has no more than one nonzero entry in any row and any column, $m(\gamma, h, D')$ equals to the number of zero rows (columns). Introduce notations
\begin{equation}\label{corank}
    m(D, h, D') = \sum_{\gamma\in D} m(\gamma, h, D'),\quad \quad s(D,D') = |D|+|D\setminus D'|.
\end{equation}
For $D=\varnothing$, we take $m(D, h, D')=0$.\\
\Theorem\Num\label{vvv}. The value of supercharacter on the superclass equals to
\begin{equation}\label{value}
\chi_{\theta,D}(K_{h,D'})  = \delta(D,h,D')(-1)^{|D\cap D'|}q^{m(D, h, D')}(q-1)^{s(D,D')}\theta(h).
\end{equation}
\Proof. Denote  by  $\phi(\theta,D,h,D')$ the expression on the right hand side of the formula (\ref{value}).  Simplify notations  $\chi=\chi_{\theta,D}$,~ $\xi=\xi_{\theta,D}$, ~ $g=g_{h,D'}$.
We calculate the induced character  $\chi(g) = \ind(\xi, G_\la,G)$ using well known formula
\begin{equation}\label{chis}
   \chi(g) = \frac{1}{|G_\la|}\sum
\xi(s^{-1}gs),
 \end{equation}
 where $ s^{-1}gs\in G_\la$. If the set  of  $\{s\}$, obeying  $ s^{-1}gs\in G_\la$,  is empty, we take  $\chi(g)=0$.

 Let  $D=\{\gamma_1,\ldots, \gamma_r\}$.
For each $1\leq k\leq r$ we denote
 $$\la_k = \sum_{a=k}^r E_{\gamma_a}^*, \quad D_k = \{\gamma_1,\ldots\gamma_k\},\quad\quad D_0=\varnothing.$$
Here $\la_{1}=\la_D$ and $D_r=D$. Set $\la_{r+1} =0$.
Notice that $\la_k$ is a sum of the elements of the dual basis  $E^*_\gamma$ over  $\gamma\in D\setminus D_{k-1}$.

The subgroups $G_k=G_{\la_k}$ form the chain
$$G_\la = G_1\subset G_2\subset\ldots\subset G_r \subset G_{r+1} = G.$$
Define  the character  $\chi_k$ of the group  $G_k$ by the formula
$$\chi_k =\ind(\xi, G_\la,G_k).$$
 Then $\chi_1 = \xi$ is a character of the subgroup  $G_\la=G_1$. For each $1\leq k\leq r$ we have
$$ \chi_{k+1} = \ind(\chi_k, G_k,G_{k+1}).$$
We shall prove
\begin{equation}\label{induction}
 \chi_{k+1}(g) = \phi(\theta,D_{k},h,D')\varepsilon^{\la_{k+1}(x_{D'})},
\end{equation}
using the induction method with respect to  $0\leq k\leq r$.
For  $k=0$, the formula  (\ref{induction}) is true, since $\chi_1(g)=\theta(h)
\varepsilon^{\la(x_{D'})}$.
Suppose that  (\ref{induction}) is true for  $k-1$; let us prove for  $k$.

Let  $\gamma=(i,j)$, ~$1\leq i<j\leq n$ be an arbitrary positive root. Let $P'(\gamma)$ be the submatrix of the matrix $g-1$ with the systems of rows
 $\{i\}\cup [i+1,j-1]$ and columns $ [i+1,j-1]$.
 Respectively, $P''(\gamma)$ is the submatrix of the matrix  $g-1$ with the systems of rows $[i+1,j-1]$ and columns $ [i+1,j-1]\cup \{j\}$.

   The induction assumption and formula  (\ref{chis}) implies
 \begin{equation}\label{qqq}
       \chi_{k+1}(g) = \phi(\theta,D_{k-1},h,D')  M(\gamma_k,h,D)  \varepsilon^{\la_{k+1}(x_{D'})},
  \end{equation}
where
\begin{equation}\label{MMM}
M(\gamma_k,h,D) = \delta_0(\gamma_k,h)\sum_{t_i,t_j\in\Fq^*}    \varepsilon^{t_i^{-1}t_j}\sum_{(t_i,\overline{x})P'(\gamma_k)=0} \varepsilon^{(\overline{x},\overline{p})},
\end{equation}
$\overline{x}=(x_1,\ldots, x_d)$, ~$\overline{p}=(p_1,\ldots,p_d)$ is the last column of the submatrix  $P''(\gamma_k)$, ~~ $(\overline{x},\overline{p})=x_1p_1+\ldots+x_dp_d$ ~and ~ $d=j-i-1$.

If there exists  $\gamma'\in D'$ lying in  $\De(\gamma_k)$ (i.e.  $\delta'(\{\gamma_k\},h,D')=0$), then  $$\rank P(\gamma_k)< \rank P'(\gamma_k).$$
The underlying subset of the second sum in  (\ref{MMM}) is empty. Therefore  $$M(\gamma_k,h,D)=0.$$

If there is no  $\gamma'\in D'$ lying in  $\De(\gamma_k)$ (i.e.  $\delta'(\{\gamma_k\},h,D')=1$), then the first row of the submatrix $P'(\gamma_k)$ is zero. Hence
\begin{equation}\label{MsMs}
M(\gamma_k,h,D) = \delta_0(\gamma_k,h)\delta'(\{\gamma_k\},h,D')\left(\sum_{t_i,t_j\in\Fq^*}    \varepsilon^{t_i^{-1}t_j}\right) \left(\sum_{\overline{x}P(\gamma_k)=0} \varepsilon^{(\overline{x},\overline{p})}\right)
\end{equation}
We calculate the third factor  $M_3$ in (\ref{MsMs}):
$$M_3 = \sum_{t_i,t_j\in\Fq^*}    \varepsilon^{t_i^{-1}t_j}
=\left\{\begin{array}{l}
(q-1)^2, \quad\mbox{if}~~ \gamma_k\notin D',\\
(q-1)(-1), ~~ \mbox{otherwise}.
\end{array}\right.$$
Then $$ M_3 =(-1)^{|\{\gamma_k\}\cap D'|}
(q-1)^{s(\{\gamma_k\}, D')}.$$
To calculate the forth factor in (\ref{MsMs}) we notice that
$$\sum_{\overline{x}\in W}\varepsilon^{(\overline{x},\overline{p})} =\left\{\begin{array}{l}
|W|, \quad\mbox{if}~~ \overline{p}\in W^\perp,\\
0, ~~ \mbox{otherwise}.
\end{array}\right.$$
for any linear subspace $W$.
We calculate the forth factor  $M_4$ in (\ref{MsMs}):
$$M_4 = \sum_{\overline{x}P(\gamma_k)=0} \varepsilon^{(\overline{x},\overline{p})} =
\left\{\begin{array}{l}
q^{\mathrm{corank} P(\gamma_k)}, \quad\mbox{if}~~ \overline{p}\in \mathrm{Im}(P(\gamma_k)) ,\\
0, ~~ \mbox{otherwise}.
\end{array}\right.$$
Then
$$ M_4 = \delta''(\{\gamma_k\},h,D')q^{m(\{\gamma_k\},h,D')}.$$
 Since $$\delta_0(\gamma_k,h)\delta'(\{\gamma_k\},h,D')\delta''(\{\gamma_k\},h,D') = \delta(\{\gamma_k\},h,D'),$$ after substitution  of   $M_3$ and $M_4$ into  (\ref{MsMs}) we obtain:
\begin{equation}\label{Mcalc}
    M(\gamma_k,h,D) = \delta(\{\gamma_k\},h,D')(-1)^{|\{\gamma_k\}\cap D'|}
(q-1)^{s(\{\gamma_k\}, D')}q^{m(\{\gamma_k\},h,D')}.
\end{equation}
Substituting  (\ref{Mcalc}) into  (\ref{qqq}), we have got
$$  \chi_{k+1}(g) = \phi(\theta,D_{k},h,D')\varepsilon^{\la_{k+1}(x_{D'})},$$
this proves  (\ref{induction}) for $k$.

  Taking $k=r$ in (\ref{induction}), we finally obtain $\chi(g) = \phi(\theta,D,h,D')$.
$\Box$


\begin{thebibliography}{99}

\bibitem{A1}
Andr\'e C.A.M., {\it Basic Sums of Coadjoint Orbits of the Unitriangular group},
Journal of Algebra {\bf 176} (1995), 959-1000.
\bibitem{A2}
Andr\'e C.A.M., {\it Basic  character table of the Unitriangular group}, Journal of Algebra {\bf 241} (2001),  437-471.
\bibitem{A3}
Andr\'e C.A.M., {\it Hecke algebra for the basic representations of the
unitriangular group}, Proc. Amer. Math. Soc. {\bf 132} (2003), no. 4,
987-996.
\bibitem{DI}
Diaconis P., Isaacs I.M., {\it Supercharacters and  superclasses for
algebra groups}, Trans.Amer.Math.Soc. {\bf 360} (2008),  2359-2392.
\bibitem{Yan}
Ning Yan, {\it Representation Theory of finite unipotent linear groups}, Dissertation, 2001 (arXiv: 1004.2674).
\bibitem{VERY}
 Aguiar M.,  Andr\`{e} C.,   Benedetti C., Bergeron N., Zhi Chen,  Diaconis P.,  Hendrickson A.,  Hsiao S.,  Isaacs I.M.,  Jedwab A.,  Johnson K., Karaali G.,  Lauve A., Tung Le,  Lewis S., Huilan Li,  Magaarg K.,  Marberg E.,   Novelli J-Ch., Amy Pang,  Saliola F.,  Tevlin L., Thibon  J-Y.,  Thiem N.,  Venkateswaran V., Vinroot C.R., Ning Yan,  Zabricki M., {\it Supercharacters, symmetric functions in noncommuting variables, and related Hopf algebras}, Advances in Mathematics {\bf 229}(2012), no.4,  2310-2337.

\bibitem{H}
Hendrickson A.O.F., {\it Supercharacter theory costructions corresponding to Schur ring products}, Comm. Algebra {\bf 40}(2012), no.12, 4420-4438.

\bibitem{SA}
 Andrews S., {\it Supercharacter theory constructed by the method of little groups}, 2014, available at arXiv 1405.5472.
\bibitem{Lang}
Lang A., {\it Supercharacter theories and semidirect products}, 2014, available at arXiv 1405.1764.

\bibitem{Pi}
Pierce R.S., {\it Associative algebras}, Springer-Verlag, New York, 1982.
\bibitem{DK}
Drozd Yu.N., Kirichenko V.V., {\it Finite dimensional algebras},  Vischa Shkola, Kiev,  1980 [in russian].
\bibitem{CR}
  Curtis Ch. W.,  Reiner I., {\it Representation theory of finite groups and associative algebras}, Interscience Publishers,  New York, 1962.
\bibitem{Serr}
 Serr J.-P., {\it Linear representations of finite groups}, New-York, Springer-Verlag, New-York, 1977.





\end{thebibliography}
\end{document}